\documentclass{article}

\usepackage{arxiv}
\usepackage{amssymb}
\usepackage{amsmath}
\usepackage{empheq}
\usepackage{amscd}
\usepackage{graphicx}
\usepackage{graphics}
\usepackage[noadjust]{cite}
\usepackage{amsthm}
\usepackage{caption}
\usepackage{multirow}
\usepackage{array}
\usepackage{rotating}

\usepackage{indentfirst}
\usepackage{tikz}
\usepackage{calc}
\usepackage{epsfig}
\usepackage{natbib}
\usepackage[makeroom]{cancel}
\usepackage{multicol}
\usepackage{csquotes}
\usepackage{algorithm}
\usepackage{algorithmic}
\usepackage{enumitem}
\usepackage{hyperref}
\usepackage{url}
\hypersetup{colorlinks,linkcolor={blue},citecolor={blue},urlcolor={blue}}

\usepackage{timet}
\usepackage{graphicx}
\usepackage{multirow}
\usepackage{multicol}
\usepackage{lipsum}
\usepackage{timet}
\usepackage{epsfig}
\usepackage{amsmath}
\usepackage{amsfonts}
\usepackage{amssymb}
\usepackage{color}
\numberwithin{equation}{section}
\usepackage{caption}
\usepackage{float}
\usepackage{subcaption}
\usepackage{graphics}
\usepackage{xcolor,graphicx}
\usepackage{csquotes}
\usepackage{mathtools}
\usepackage{optidef}

\newcommand{\vertiii}[1]{{\left\vert\kern-0.25ex\left\vert\kern-0.25ex\left\vert #1 
    \right\vert\kern-0.25ex\right\vert\kern-0.25ex\right\vert}}

\usepackage{hyperref}       

\renewcommand{\u}[0]{\bold{u}}                        

\renewcommand{\b}[0]{\bold{b}}                        
\renewcommand{\d}[0]{\bold{d}}                        
\renewcommand{\P}[0]{\bold{P}}                        

\newcommand{\m}[0]{\bold{m}}                        
\newcommand{\x}[0]{\bold{x}}                        
\newcommand{\w}[0]{\omega}                          
\newcommand{\A}[0]{\bold{A}}                        
\newcommand{\del}[0]{\bold{\nabla}^2}               
\newcommand{\Diag}[0]{\text{Diag}}                  
\newcommand{\AL}[0]{\mathcal{L}}                    
\newcommand{\penaltyparb}[0]{\lambda}               
\newcommand{\penaltypard}[0]{\mu}                   
\newcommand{\dualb}[0]{\bold{v}}                    
\newcommand{\duald}[0]{\bold{w}}                    
\newcommand{\reg}[0]{\mathcal{R}}                   
\newcommand{\Nr}[0]{r}                         
\newcommand{\N}[0]{n}                          
\newcommand{\Loc}[0]{\boldsymbol{\Phi}}                   
\newcommand{\s}[0]{\bold{s}}                        
\newcommand{\spoints}[0]{p}                        
\newcommand{\fpoints}[0]{q}                        
\renewcommand{\H}[0]{\bold{H}}                   
\newcommand{\grad}[0]{\bold{g}}                   

\newcommand{\eref}[1]{Eq. \eqref{#1}}  

\begin{document}
\title{ADMM-based full-waveform inversion for microseismic imaging}

\author{\href{http://orcid.org/0000-0003-1805-1132}{\includegraphics[scale=0.06]{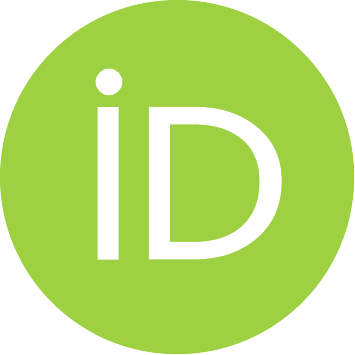}\hspace{1mm}Hossein S. Aghamiry} \\
  University Cote d'Azur - CNRS - IRD - OCA, Geoazur, Valbonne, France. 
  \texttt{aghamiry@geoazur.unice.fr}
\And
  \href{https://orcid.org/0000-0002-9879-2944}{\includegraphics[scale=0.06]{orcid.pdf}\hspace{1mm}Ali Gholami} \\
  Institute of Geophysics, University of Tehran, Tehran, Iran.
  \texttt{agholami@ut.ac.ir} \\ 
  \And
\href{http://orcid.org/0000-0002-4981-4967}{\includegraphics[scale=0.06]{orcid.pdf}\hspace{1mm}St\'ephane Operto} \\ 
  University Cote d'Azur - CNRS - IRD - OCA, Geoazur, Valbonne, France. 
  \texttt{operto@geoazur.unice.fr}
\And
\href{http://orcid.org/0000-0003-2775-7647}{\includegraphics[scale=0.06]{orcid.pdf}\hspace{1mm}Alison Malcolm} \\ 
  Department of Earth Sciences, Memorial University of Newfoundland, St John’s NL A1B 3X5, Canada. 
  \texttt{amalcolm@mun.ca} 
  }

\renewcommand{\shorttitle}{Microseismic imaging using IR-WRI, Aghamiry et al.}

\maketitle

\begin{abstract}
Full waveform inversion (FWI) is beginning to be used to characterize weak seismic events at different scales, an example of which is microseismic event (MSE) characterization. However, FWI with unknown sources is a severely underdetermined optimization problem, and hence requires strong prior information about the sources and/or the velocity model. The frequency-domain wavefield reconstruction inversion method (WRI) has shown promising results to mitigate the nonlinearity of the FWI objective function that is generated by cycle-skipping. WRI relies on the reconstruction of data-assimilated wavefields, which approach the true wavefields near the receivers, a helpful feature when the source is added as an additional optimization variable.  
We present an adaptation of a recently proposed version of WRI based on the alternating direction method of multipliers (ADMM) that first finds the location of the MSEs and then reconstructs the wavefields and the source signatures jointly. Finally, the subsurface model is updated to focus the MSEs at their true location. The method does not require prior knowledge of the number of MSEs. The inversion is stabilized by sparsifying regularizations separately tailored to the source location and velocity model subproblems. The method is tested on the Marmousi model using one MSE and two clusters of MSEs with two different initial velocity models, an accurate one and a rough one, as well as with added noise. In all cases, the method accurately locates the MSEs and recovers their source signatures. 
\end{abstract}

\section{INTRODUCTION}
Full waveform inversion (FWI) is a waveform matching procedure that is generally used to reconstruct subsurface models with wavelength resolution \citep{Virieux_2009_OFW}. 
The FWI optimization problem is commonly solved assuming that the positions of the controlled sources are known, while the temporal source signature is estimated together with the subsurface model by variable projection or in an alternating scheme \citep{Pratt_1999_SWIb,Plessix_2011_GJI,Aravkin_2012_SEF,Rickett_2013_VPM}.\\
This paper focusses on the case where the source location and time signature are not known.  The source can be active (e.g., explosion) or passive (e.g., earthquake, microseismic). Our numerical examples focus on characterizing weak seismic events at different scales, especially for microseismic event (MSE) imaging, but the method can easily be tailored to tectonic earthquakes.  When estimating the source, we decompose it into a temporal signature and a spatial location function, both of which strongly influence the FWI results.  \citet{Lee_2003_SIF} show that a small source error can seriously affect the estimated model parameters, and so source estimation is an important step in successful FWI.

During oil and gas production, CO$_2$ injection, and geothermal applications, fluid injections generate fracturing and cause small earthquakes (microseismicity). Locating the MSEs and studying their focal mechanisms are important in monitoring flow mobility and optimizing production.  The seismic waves generated by these earthquakes are recorded continuously at the surface or in wells, and their traveltimes or waveforms are inverted to estimate the source time function (source signature) and the location of the MSEs, assuming a known subsurface model. This setting implies that the recorded wavefield has been triggered by a blended source with contributions from of all of the MSEs.

Pioneering methods for microseismic location rely on picking the traveltimes of different P and S components to locate the MSEs \citep{Thurber_2000_AIS,Han_2009_HLU}.  These methods are time-consuming and may not be robust in the presence of noise \citep{Thurber_2000_AIS,Han_2009_HLU}. 
More recently, wavefield-based imaging techniques have received a lot of attention for the relocation of earthquakes and MSEs \citep{McMechan_1982_DSP,Gajewski_2005_RMS,Michel_2014_GCF, Sjogreen_2014_SEF,Kaderli_2015_MEE,DeRidder_2018_FWI,Shanan_2019_FSM,Shekar_2019_FIM,Song_2019_MES}. The most basic approaches rely on time-reverse imaging, where the data are propagated backward in time using the adjoint of the wave-equation operator \citep{McMechan_1982_DSP,Gajewski_2005_RMS,Nakata_2016_RTM,Shanan_2019_FSM}. The source location and the origin time are found by tracking where, in both space and time, the maximum focusing occurs. These approaches suffer from a lack of resolution and, as in any migration technique, require a kinematically accurate background velocity model.  To mitigate these problems, there has been some interest in applying FWI for imaging such microseismicity \citep{Montgomery_2010_HFH,Li_2020_RAC}. 
Moving from migration-based techniques to FWI-based techniques opens the door to estimate the source signatures, the source locations and in some cases the subsurface properties, by solving a multivariate optimization problem.
This is possible because all of the wave information is involved in the FWI procedure. This multivariate optimization problem attempts to reconstruct wavefields, sources, and model parameters from partial measurements of the wavefield. 
However, this problem is highly under-determined, even when the subsurface parameters are processed as inactive parameters (i.e. are not updated). To mitigate this, the method requires some priors like the number of events, their approximate locations, the sparsity of the source distribution or an accurate background velocity model to tighten the null space of the inverse problem \citep{Kaderli_2015_MEE,Kaderli_2018_SVF,Shekar_2019_FIM, Shanan_2019_FSM}.
\citet{Michel_2014_GCF} reviewed how to compute the FWI gradient for source parameters in VTI media. 
\citet{Kaderli_2015_MEE,Shekar_2019_FIM, Shanan_2019_FSM} proposed robust FWI algorithms with sparsifying regularization to refine the location of MSEs and estimate their temporal signatures. Their approaches first estimate the location of the events by wavefield extrapolation with sparsifying regularization. Following this, the source signatures of the identified events are reconstructed separately.  These approaches do not require assumptions about the number or the nature of the events, but require a fairly accurate background velocity model. \\
A few studies have been presented to jointly locate the events and update the subsurface velocity model \citep{Sun_2016_FPF,Song_2019_MES}. \citet{Sun_2016_FPF} update the source and the velocity model in an alternating manner. The descent direction of the under-determined source estimation problem is preconditioned by a weighting of the source term inferred from a cross-correlation time-reversal imaging condition. \citet{Song_2019_MES} build the source image using the imaging condition of \citet{Nakata_2016_RTM} and then update the velocity model by penalizing the energy of the source image away from the estimated source location with an annihilator. This approach requires one to process each event separately and hence to separately identify each event in the data. 
\citet{Kaderli_2015_MEE} estimate the MSE locations and source signatures in an alternating way when the source is considered as a product of independent temporal source function and spatial source location function, while \citet{Shanan_2019_FSM} update both of them jointly.  \\
%
%
We extend these past works by using a more general form of FWI. In this more general form, FWI can be cast as a constrained optimization problem that aims to estimate the wavefields and the subsurface parameters by fitting the recorded data subject to the constraint that the wave equation is satisfied \citep{Haber_2000_OTS}. This approach was developed because even when the source location and time signature are known, it is well established that FWI is highly nonlinear. Part of this nonlinearity can be viewed as arising when the full search space encompassed by the wavefields and the subsurface parameters is projected onto the subsurface parameter space. This happens via an elimination of the wavefield variables, by assuming that the wavefields exactly satisfy the wave-equation at each FWI iteration. This variable elimination makes FWI prone to cycle skipping when the initial model is not accurate enough to predict recorded traveltimes with an error smaller than half a period \citep{Virieux_2009_OFW}. To avoid this projection, some approaches implement the wave equation as a soft constraint with a penalty method such that the data can be closely matched with inaccurate subsurface models from the early FWI iterations by not requiring that the wave equation be satisfied exactly \citep{Abubakar_2009_FDC,VanLeeuwen_2013_MLM,vanLeeuwen_2016_PMP}. Then, the subsurface model is updated by solving an overdetermined quadratic optimization problem, which consists of minimizing the source residuals generated by the relaxation of the constraint that the wavefields exactly solve the wave equation. In these extended approaches, the wavefields are reconstructed by solving,  in a least-squares sense, an overdetermined linear system comprised of the wave equation weighted by the penalty parameter and the observation equation relating the simulated wavefield to the data through a sampling operator. In other words, the wavefields are reconstructed with data assimilation, which makes them approach the true wavefields near the receivers, a helpful feature when the source is added as a new variable. This approach was called Wavefield Reconstruction Inversion (WRI) by \citet{VanLeeuwen_2013_MLM}. A variant of WRI, based upon the method of multipliers or augmented Lagrangian method,  was proposed by \citet{Aghamiry_2019_IWR} to increase the convergence rate and decrease the sensitivity of the algorithm to the choice of the relaxation (penalty) parameter. The augmented Lagrangian method combines a penalty method and a Lagrangian method, where the penalty term is used to implement the initial relaxation of the constraint, and the Lagrangian term automatically tunes the sensitivity of the optimization to the constraint in iterations. The Lagrange multipliers are updated with gradient ascent, which controls the constraint violations. This method is called Iteratively-Refined(IR)-WRI, where the prefix IR refers to the iterative refinement ( i.e. defect correction) action of the Lagrange multipliers. \\
%
%
In this study, we propose a new microseismic imaging algorithm for event location and velocity model building based on IR-WRI \citep{Aghamiry_2019_IWR}. 
We extend the IR-WRI method to solve for the signatures and locations of the MSEs as additional variables. 
Beginning from the initial velocity model, without any assumptions about the MSEs, the first data assimilated wavefield is reconstructed for a band of frequencies starting without a source term (namely, only the data drive the wavefield reconstruction). 
Then by using the extracted multi-frequency data-assimilated wavefields, we estimate a mean source term averaged over frequencies. 
During this mean source estimation, we use sparsifying denoising to focus the blended source and hence further facilitate the localization of the MSEs.
After a few iterations of this two-step process alternating between wavefield reconstruction and mean source estimation, we apply a peak finder algorithm to the final predicted mean source to extract the location of the MSEs. Then, we jointly update the data assimilated wavefields as well as the source signatures of the picked MSEs keeping the velocity model as a fixed parameter. Finally, we update the velocity model by minimizing the wave-equation errors when the wavefields as well as the locations and signatures of MSEs are kept fixed.
The proposed algorithm does not require assumptions about the number or type of MSEs and their locations. However, it does require that the velocity model contains the low wavenumber components of the model. \\
We first review the different steps of the method. Then, we illustrate the method with the Marmousi synthetic example. Starting from an accurate version of the Marmousi model, we show how the method manages to locate a single source without updating the velocity model. Then, we repeat the same test with a highly-smoothed starting velocity model and show how velocity model updating allows for the accurate location of the event. Then, we complicate the latter test with two small clusters of point sources and the results confirm the potential of IR-WRI for MSE localization.  As a last example, we test the method on data with added random noise.
%
\section{METHOD}
This paper relies on the frequency-domain formulation of FWI. Accordingly, we review the method with a discrete matrix formalism \citep{Pratt_1998_GNF}. \\
In passive experiments, the source $\b$ is unknown and can be approximated as a superimposition of $\spoints$ point sources (resembling a blended source).
Accordingly, the source vector $\b(\w)$ for frequency $\w$ reads
\begin{equation} \label{source_term}
\b(\w)=\sum_{j=1}^{\spoints} [\s(\w)]_j \delta(\x-\x_j)=\Loc\bold{s}(\w),
\end{equation} 
where $\delta(\x)$ is the delta function, $\x_j$ is the point source position, $[\s(\w)]_j \in \mathbb{C}^{\spoints \times 1}$ denotes the source signature, at angular frequency $\w$, associated to the $j$th point source (at location $\x_j$), $\Loc \in \mathbb{R}^{\N \times \spoints}$ is a tall matrix, the columns of which contain shifted delta functions at the positions of the MSEs ($\N$ is the number of discretization points in the model). The goal of microseismic imaging is to find the MSE location matrix $\boldsymbol{\Phi}$, the MSE signature vector $\bold{{s}}_\omega$ and the velocity model (provided that the data set provides a sufficient illumination of the model). \\
Frequency-domain FWI with an unknown blended source $\bold{b}(\omega)$ can be written as 
\begin{mini} 
{\substack{\m\in \mathcal{M},\\ \u(\w_1),\b(\w_1),\hdots,\\ \u(\w_q),\b(\w_q)}}{\reg_m(\m)+\sum_{\w=w_1}^{\w_q}\reg_b(\b(\w))}
{\label{init}}{}
\addConstraint {\A(\m,\w)\u(\w)}{=\b(\w), \quad}{\w =\w_1 ,\ldots ,\w_\fpoints}
\addConstraint {\P\u(\w)}{=\d(\w),\quad}{\w =\w_1 ,\ldots ,\w_\fpoints}
\end{mini} 
where $\m$ is the squared slowness, $\reg_m$ and $\reg_b$ are appropriate regularization functions, 
 $\A(\m,\w)=\del+\w^2 \Diag(\m) \in \mathbb{C}^{\N\times \N}$ is the Helmholtz operator, $\del$ is the Laplacian operator, $\Diag(\cdot)$ denotes a diagonal matrix with $\cdot$ on its main diagonal, $\u(\w) \in \mathbb{C}^{\N\times 1}$ and $\d(\w) \in \mathbb{C}^{\Nr\times 1}$ denote the wavefield and the recorded data for frequency $\w$,  respectively, $\P \in \mathbb{R}^{\Nr\times \N}$ is the observation operator that samples $\u(\w)$ at receiver locations, and $\Nr$ is the number of receivers. Finally, $\m\in \mathcal{M}$ is a bounding constraint on the model parameters where 
\begin{equation}
\mathcal{M} = \{\bold{m} \vert \bold{m}_{min} \leq \bold{m} \leq \bold{m}_{max}\}.
\end{equation}\\ 
The classical implementation of the FWI as formulated in \eqref{init} would enforce the wave-equation constraint, $\u(\w)=\A(\m,\w)^{-1}\b(\w)$, in the observation-equation constraint and process the latter as a penalty term leading to the following optimization problem:
\begin{mini} 
{\substack{ \m\in \mathcal{M},\b(\w_1),\hdots,\b(\w_q)}}{\reg_m(\m)+\sum_{\w=w_1}^{\w_q}\reg_b(\b(\w))+\frac{\penaltyparb}{2} \sum_{\w=w_1}^{\w_q}\|\bold{G}(\m,\w)\b(\w)-\d(\w)\|_2^2,}
{\label{pratt_FWI}}{}
\end{mini}
%
where $\penaltyparb$ is the penalty parameter and $\bold{G}(\m,\w)=\P\A(\m,\w)^{-1}$.  
The problem in equation \eqref{pratt_FWI} is severely underdetermined due to the unknown source and highly nonlinear due to the oscillatory nature of the Green's functions, which makes the waveform inversion prone to cycle skipping. 
\citet{Michel_2014_GCF} solved this problem when they assume a good initial estimate of the number of the point sources, $p$, and their approximate locations, $\Loc$.  \citet{Kaderli_2015_MEE,Shekar_2019_FIM, Shanan_2019_FSM}  used the sparsity promoting $\ell_1$-norm regularization, $\reg_b(\b(\w))=\|\b(\w)\|_1$, to enforce the sparsity of the source term, i.e., predicting the data with a minimum number of point sources. They solved the problem using the time domain formulation of FWI.

In this paper, we extend the iteratively-refined wavefield reconstruction inversion (IR-WRI) \citep{Aghamiry_2019_IWR} to solve problem \eqref{init}. IR-WRI relies on the augmented Lagrangian method, which combines the penalty method with the Lagrangian method \citep[][ Chapter 17]{Nocedal_2006_NO}. 
%
The augmented Lagrangian function associated with problem in \eref{init} is given by
\begin{align} \label{eqpsi}
\AL(\m,\{\u(\w)\},\{\b(\w)\},\{\dualb(\w)\},\{\duald(\w)\}) &=
\reg_m(\m)+\sum_{\w=w_1}^{\w_q}\reg_b(\b(\w))  \\
&+ \sum_{\w=w_1}^{\w_q}\langle \dualb(\w),\A(\m,\w)\u(\w)-\b(\w)\rangle  
+ \sum_{\w=w_1}^{\w_q}\langle \duald(\w),\P\u(\w)-\d(\w)\rangle  \nonumber \\
& + \frac{\penaltyparb}{2} \sum_{\w=w_1}^{\w_q}\|\A(\m,\w)\u(\w)-\b(\w)\|_2^2 + 
\frac{\penaltypard}{2}\sum_{\w=w_1}^{\w_q}\|\P\u(\w)-\d(\w)\|_2^2, \nonumber 
\end{align}
where the scalars $\penaltyparb,\penaltypard>0$ are the penalty parameters assigned to the wave equation and the observation equation constraints, respectively, and $\dualb(\w)$ and $\duald(\w)$ are the Lagrange multipliers.
%
Beginning with an initial model $\m^{0}$, $\b^0(\w)=0~\forall \w$, we compute an initial set of monochromatic wavefields $\u^{0}(\w)$ by solving the following overdetermined systems in a least-squares sense:
\begin{equation} \label{close_U0}
\left(
\begin{array}{c}
\sqrt{\penaltyparb} \A(\m^0,\omega) \\
\sqrt{\penaltypard}\P \\
\end{array}
\right)
\bold{u}^{0}(\omega)
=
\left(
\begin{array}{c}
0 \\
\sqrt{\penaltypard}\d(\omega)
\end{array}
\right), ~~~~~~~~~~~~~~~~~~~~ \omega=\omega_1~,....,~\omega_q.
 \end{equation} 
Then, beginning with $\dualb^{0}(\omega)=\bold{0}$ and $\duald^{0}(\omega)=\bold{0}, \forall \omega$, 
we solve the multivariate optimization problem, equation \eqref{eqpsi}, iteratively by using the ADMM as
\begin{subequations}
\label{ADMM}
 \begin{empheq}[left={\empheqlbrace\,}]{align}
(\{\u(\w)^{k+1}\},\{\b(\w)^{k+1}\})=& \underset{\{\u(\w)\},\{\b(\w)\}}{\arg\min} ~ \AL(\m^k,\{\u(\w)\},\{\b(\w)\},\{\dualb(\w)^k\},\{\duald(\w)^k\}) \label{primal_sig+wavefield}\\
\bold{m}^{k+1} &= \underset{\bold{m}\in \mathcal{M}}{\arg\min} ~ \AL(\m,\{\u(\w)^{k+1}\},\{\b(\w)^{k+1}\},\{\dualb(\w)^k\},\{\duald(\w)^k\}) \label{primal_sigma}\\
\dualb(\omega)^{k+1} &= \dualb(\omega)^{k}  + \penaltyparb (\A(\m^{k+1},\w)\u(\omega)^{k+1} - \b(\omega)^{k+1}), ~\omega=\omega_1~,....,~\omega_q \label{dual_b}\\ 
\duald(\omega)^{k+1} &= \duald(\omega)^{k}   + \penaltypard (\P\u(\omega)^{k+1} - \d(\omega)), ~\omega=\omega_1~,....,~\omega_q  \label{dual_d}
\end{empheq}
\end{subequations} 
where $k$ is the (outer) iteration number. The penalty parameters $\penaltyparb,\penaltypard>0$ are tuned such that a dominant weight $\penaltypard$ is given to the observation equation at the expense of violating the wave equation during the early iterations to guarantee a data fit that prevents cycle skipping at receivers even at early iterations. The iterative update of the Lagrange multipliers progressively corrects the errors introduced by these penalizations such that both the observation equation and the wave equation are satisfied at the convergence point with acceptable accuracies.
In the next two subsections, we review how to solve each of optimization subproblems \eqref{primal_sig+wavefield}-\eqref{primal_sigma}.
\subsection{Estimation of the sources and wavefields}
\begin{figure}
\centering
\includegraphics[width=0.7\textwidth]{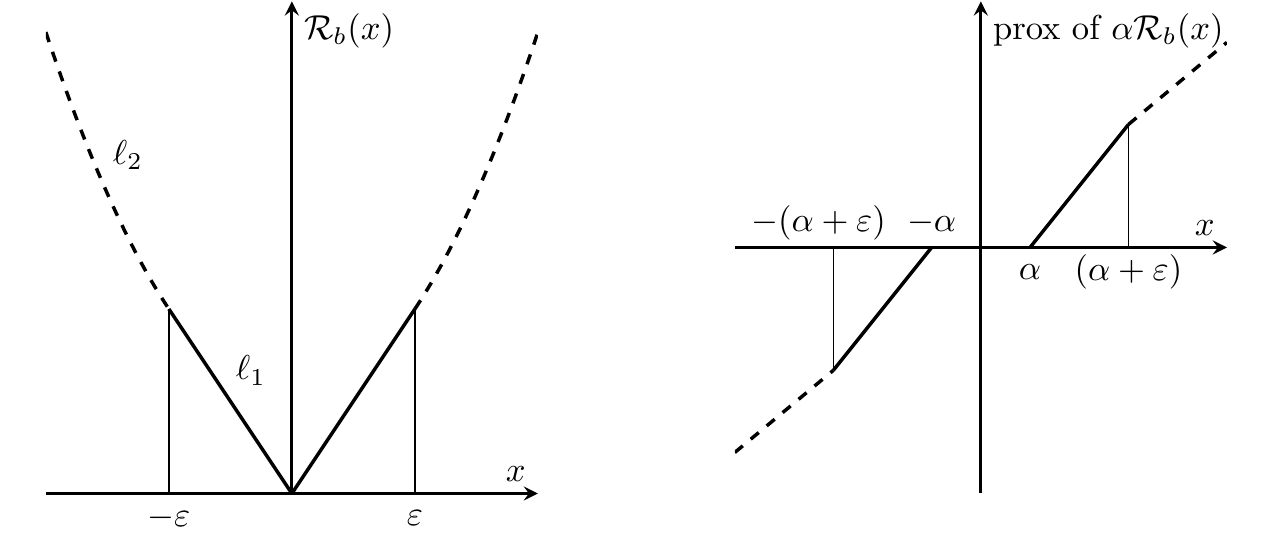}
\caption{(left) Berhu regularization function and (right) the associated proximity operator. The dashed lines are due to the $\ell_2$ norm, the solid lines are due to the $\ell_1$ norm.}
\label{fig:Berhu}
\end{figure}

Due to the ill-conditioning of the problem of estimating $\{\u(\w)^{k+1}\}$ and $\{\b(\w)^{k+1}\}$, we solve the subproblems \eqref{primal_sig+wavefield} in a two step manner. 
Since the event location matrix $\Loc$ in Eq. \eqref{source_term} is frequency independent, we first solve the optimization problem for a mean source $\bar{\b}^k=\frac{1}{q}\sum_{\w=w_1}^{\w_q}\b^k(\omega)$ averaged over the frequency to reduce the search space.
 The optimization problem over $\bar{\b}$ reads 
\begin{align} \label{meansource}
\bar{\b}^{k+1}&=\underset{\substack{\bar{\b}}}{\arg\min} ~~
\reg_{b}(\bar{\b})
+ \sum_{\w=w_1}^{\w_q}\langle \dualb(\w)^{k},\A(\m^k,\w)\u(\w)^{k}-\bar{\b}\rangle  
+ \frac{\penaltyparb}{2} \sum_{\w=w_1}^{\w_q}\|\A(\m^{k},\w)\u(\w)^{k}-\bar{\b}\|_2^2.
\end{align}
By adding and subtracting the term $\|\dualb(\w)^{k}\|^2_2$ to problem \eqref{meansource}, we have \citep[, Appendix A]{Aghamiry_2019_IBC}
\begin{align}
\bar{\b}^{k+1}&=\underset{\substack{\bar{\b}}}{\arg\min} ~~
\reg_{b}(\bar{\b})
+ \frac{\penaltyparb}{2} \sum_{\w=w_1}^{\w_q}\|\A(\m^{k},\w)\u(\w)^{k}-\bar{\b}+\frac{1}{\penaltyparb}\dualb(\w)^{k}\|_2^2 \label{obj_phi}\, ,
\end{align}
where we have ignored the $-\|\dualb(\w)^{k}\|^2_2$-term as it does not impact the optimization result.
%
%
Equation \ref{obj_phi} is a denoising/proximity problem \citep{Parikh_2013_PA} applied to $\bar{\b}$, i.e.
\begin{equation} \label{measo}
\bar{\b}^{k+1}=\underset{\substack{\bar{\b}}}{\arg\min} ~~
\reg_{b}(\bar{\b}) + 
\frac{\penaltyparb q}{2} \|(\frac{1}{q}\sum_{\omega=\omega_1}^{\omega_q}\A(\m^{k},\w)\u(\w)^{k}+\frac{1}{\penaltyparb}\dualb(\w)^{k}) - \bar{\b}\|_2^2.
\end{equation}
%
%
%
%
We use the so called Berhu regularizer \citep{Owen_2007_ARH}, which is a hybrid function combining the $\ell_1$ norm for small values and the $\ell_2$ norm for large values:
\begin{equation}
\reg_{b}(x)=\begin{cases}
|x|~~~~~~~~~~~~ |x| \leq \varepsilon,     \\
\frac{x^2+\varepsilon^2}{2\varepsilon}  ~~~~~~ |x| > \varepsilon,
\end{cases}
\end{equation}
where $\varepsilon>0$ determines where the transition from $\ell_1$ to $\ell_2$ occurs. 
This regularizer shrinks small coefficients with the $\ell_1$ norm to promote sparsity while damping the large coefficients (Fig. \ref{fig:Berhu}). 
%
Using the Berhu regularizer, the proximity operator in Eq. \eqref{measo} admits an explicit solution for each entry of the source vector 
\begin{equation} \label{updated_source}
\bar{\b}^{k+1} =\text{prox}_{\frac{1}{\lambda q}\mathcal{B}}\left(\frac{1}{q}\sum_{\omega=\omega_1}^{\omega_q}\A(\m^{k},\w)\u(\w)^{k}+\frac{1}{\penaltyparb}\dualb(\w)^{k}\right), 
\end{equation}
where $\text{prox}_{\tau\mathcal{B}}$ is the proximity operator of the Berhu function defined as 
\begin{equation}
\text{prox}_{\alpha\mathcal{B}}(x)=\begin{cases}
\max(1-\frac{\alpha}{|x|},0)x ~~~~~~~~~ |x| \leq \alpha+\varepsilon, \\
\frac{\varepsilon}{\alpha + \varepsilon}x ~~~~~~~~~~~~~~~~~~~~~~~~~~~~~ |x|> \alpha +\varepsilon.
\end{cases}
\end{equation}
In order to improve the accuracy of the locations, we repeat the process of estimating the locations and updating the wavefields several times.
The monochromatic wavefields are updated as follows
\begin{equation} \label{close_UK}
\left(
\begin{array}{c}
\sqrt{\penaltyparb} \A(\m^{k},\w) \\
\sqrt{\penaltypard}\bold{P} \\
\end{array}
\right)
\u(\w)^{k}
=
\left(
\begin{array}{c}
\sqrt{\penaltyparb}\bar{\b}^{k+1} \\
\sqrt{\penaltypard}\d(\w)
\end{array}
\right),   ~~~~~~~~~~~~~~~~~~ \omega=\omega_1~,....,~\omega_q,
\end{equation} 
and we alternate between solving Eq. \eqref{close_UK} for $\u(\w)^{k}$ given $\bar{\b}^{k+1}$ and updating $\bar{\b}^{k+1}$ via Eq. \eqref{updated_source} using $\u(\w)^{k}$ obtained from Eq. \eqref{close_UK}. 

The final $\bar{\b}^{k+1}$ is then used to construct the location matrix $\Loc^{k+1}$. To do so, we apply a peak finder algorithm on $|\bar{\b}^{k+1}|$ to determine the location of the peaks. The number of peaks found (i.e., the number of MSE) determines the number of columns of $\Loc^{k+1}$ and the location of each peak determines the location of the delta function in the corresponding column.
Once $\Loc^{k+1}$ has been determined, the monochromatic wavefields and the time function of each MSE are updated simultaneously:
%
\begin{equation} \label{closed_u_2}
\begin{pmatrix}
\sqrt{\penaltyparb} \bold{A}_\omega^k & -\sqrt{\lambda}\boldsymbol{\Phi}^{k+1}  \\
\sqrt{\penaltypard} \bold{P} & \bold{0} 
\end{pmatrix}
\begin{pmatrix}
\u(\w)^{k+1}\\
\s(\w)^{k+1}
\end{pmatrix}
=
\begin{pmatrix}
 \frac{1}{\sqrt{\penaltyparb}}\dualb(\w)^k \\
\sqrt{\penaltypard}\d(\w)+\frac{1}{\sqrt{\penaltypard}}\duald(\w)^k
\end{pmatrix}.
\end{equation} 
%
%
\subsection{Velocity model update}
If we have adequate illumination, we can update the velocity model to improve the source locations; otherwise, this step can be skipped.
Due to the blended nature of microseismic data, the velocity update will contain noise from crosstalk. This high-frequency noise can be decreased by applying appropriate regularization.       
In this paper, we use the first-order isotropic TV regularization \citep{Rudin_1992_NTV} for $\bold{m}$. However, other regularizations such as compound or adaptive regularizations can be used in a similar way \citep[see][]{Aghamiry_2019_CRO,Aghamiry_2020_FWI}. The isotropic TV regularizer is defined as 
\begin{equation}
\mathcal{R}_m(\m)=\sum  \sqrt{(\nabla_{\!x} \m)^2 + (\nabla_{\!z}\m)^2},
\end{equation}
where $\nabla_{\!x}$ and $\nabla_{\!z}$ are respectively first-order difference operators in the horizontal and vertical directions with appropriate boundary conditions. 
The problem \ref{primal_sigma} with non-smooth TV regularization and bound constraints can be written as
\begin{equation} \label{M_sub}
\bold{m}^{k+1}= \underset{\bold{m}\in \mathcal{M}}{\arg\min} ~ \mathcal{R}_m(\m)+ \lambda
(\m^T\H_k\m - \grad_k^T \m),
\end{equation}
where
\begin{eqnarray}
&&\H_k = \sum_{\omega=\omega_1}^{\omega_q}\left(\frac{\partial \A(\m)}{\partial \m} \u(\w)^{k+1}\right)^T \left( \frac{\partial \A (\m)}{\partial \m} \u(\w)^{k+1} \right),\\
&&\grad_k = \sum_{\omega=\omega_1}^{\omega_q}\left(\frac{\partial \A (\m)}{\partial \m}\u(\w)^{k+1}\right)^T 
(\Loc^{k+1}\s(\w)^{k+1}-\frac{1}{\penaltyparb}\dualb(\w)^{k}-\del \u(\w)^{k+1}).
\end{eqnarray}
The model subproblem in Eq. \eqref{M_sub} requires a quadratic optimization with bound constrained TV regularization. There are many well documented algorithm to carry out this step \citep{Goldstein_2009_SBM,Maharramov_2015_TVM,Aghamiry_2019_IBC,Aghamiry_2019_CRO,Gholami_2019_3DD}
.\\
The proposed microseismic imaging algorithm is summarized in Algorithm~\ref{Alg2cont0}. 
It should be noted that the total algorithm consists of three main steps: (1) estimation of the source location (lines 5-9),
(2) estimation of the source signatures and the wavefields (line 10) and (3) update of the model parameters (line 11).
Step 3 can be skipped if a new sets of MSEs is processed with a sufficiently-accurate background subsurface model estimated during a prior inversion. The linear systems for data-assimilated wavefield reconstruction, equations~\ref{close_U0}, \ref{close_UK}, and \ref{closed_u_2}, have one right-hand side. Therefore, they may be solved more efficiently with preconditioned iterative solvers instead of direct solvers to tackle large computational domain. However, if multiple datasets are processed over time with the same background subsurface model, it may be beneficial to use direct solvers to compute the Lower-Upper (LU) factors once and for all once a good background model has been estimated, store them on disk and re-use them to process a new dataset efficiently by forward/backward elimination.
\begin{algorithm}[htb!]
\caption{
ADMM-based FWI for microseismic imaging.}
\label{Alg2cont0}
\scriptsize
{\fontsize{8}{8}\selectfont
\begin{algorithmic}[1]
\STATE Begin with $k=0$ and an initial model $\bold{m}^0$.
\STATE Set to zero the values of $\dualb(\w)^0$ and $\duald(\w)^0, \forall \w$.
\STATE Calculate $\bold{u}(\w)^{0}, \forall \w$ (Eq. \eqref{close_U0}). 
\item[]
\WHILE {convergence criteria  not satisfied}
\item[]
\FOR {$l=0:n_l$}
\item[]
\STATE Compute $\bar{\b}^{k+1}$ (Eq. \eqref{updated_source})
\item[]
\STATE Update $\u(\w)^{k}$, $\w=\w_1,\cdots,\w_q$ (Eq. \eqref{close_UK})
\hspace*{8.8em}%
        \rlap{\smash{$\left.\begin{array}{@{}c@{}}\\{}\\{}\\{}\\{}\\{}\\{}\\{}\end{array}\color{black}\right\}%
          \color{black}\begin{tabular}{l}Source location estimation.\end{tabular}$}}
\item[]
\ENDFOR
\item[]
\STATE Build $\boldsymbol{\Phi}^{k+1}$ from $\bar{\b}^{k+1}$ using a peak finder algorithm.
\item[]
\STATE Update $\s(\w)^{k+1}$ and $\u(\w)^{k+1}$ for $\w=\w_1,\cdots,\w_q$ (Eq. \eqref{closed_u_2}). \hspace*{1.7em}%
        \rlap{\smash{$\left.\begin{array}{@{}c@{}}\end{array}\color{black}\right\}%
          \color{black}\begin{tabular}{l}Joint update of wavefields and source signatures.\end{tabular}$}} 
\item[]
\STATE Update $\m^{k+1}$ (Eq. \eqref{M_sub})
\hspace*{17.5em}%
        \rlap{\smash{$\left.\begin{array}{@{}c@{}}\end{array}\color{black}\right\}%
          \color{black}\begin{tabular}{l}Update velocity model.\end{tabular}$}} 
\item[]
\STATE Update the dual vectors $\dualb$ and $\duald$ (Eqs. \eqref{dual_b} and \eqref{dual_d})
\hspace*{5.5em}%
        \rlap{\smash{$\left.\begin{array}{@{}c@{}}\end{array}\color{black}\right\}%
          \color{black}\begin{tabular}{l}Iterative refinement.\end{tabular}$}} 
\item[]
\STATE $k = k+1$ ,
\item[]
\ENDWHILE
 \end{algorithmic}
}
\end{algorithm}

\section{NUMERICAL EXAMPLES}
\begin{figure}
\centering
\includegraphics[width=1\textwidth]{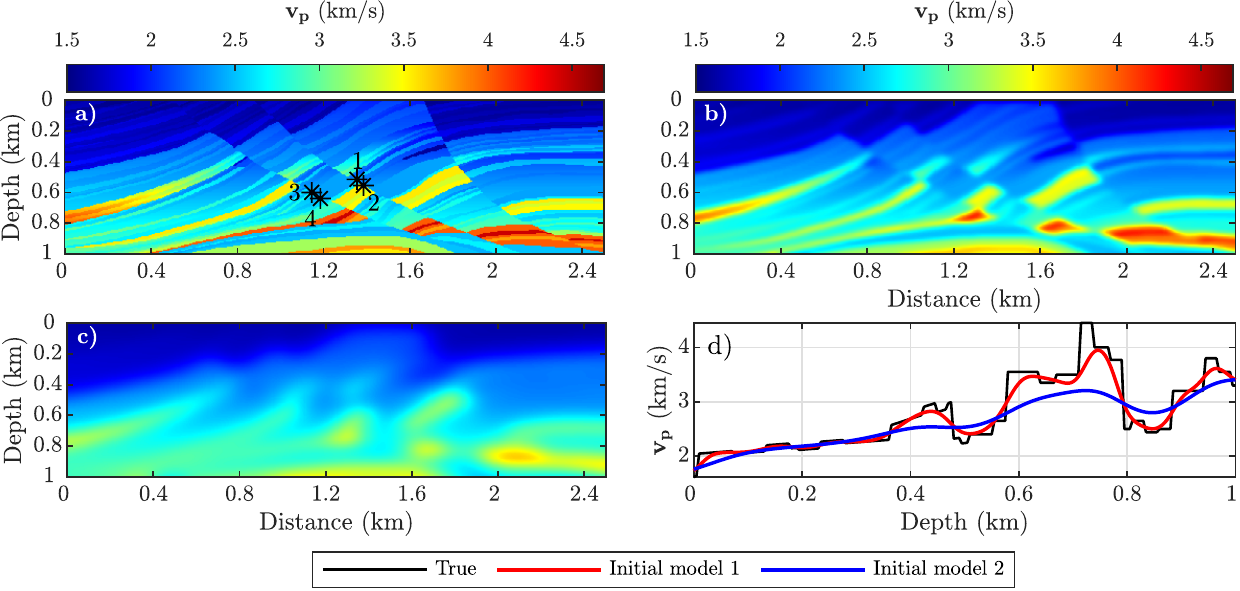}
\caption{True and initial velocity models and the locations of MSEs that are used for the numerical tests. (a) The Marmousi II velocity model, used to generate the data. The location of different MSEs, referred as MSEs 1 to 4, are shown by black asterisks. (b) Accurate initial model, called initial model 1. (c) Kinematically accurate velocity model, called initial model 2. (d) A direct comparison between true (black), initial model 1 (red), and initial model 2 (blue) at $X=1250~m$.}
\label{fig0}
\end{figure}
%
%
\begin{figure}[htb]
\centering
\includegraphics[width=1\textwidth]{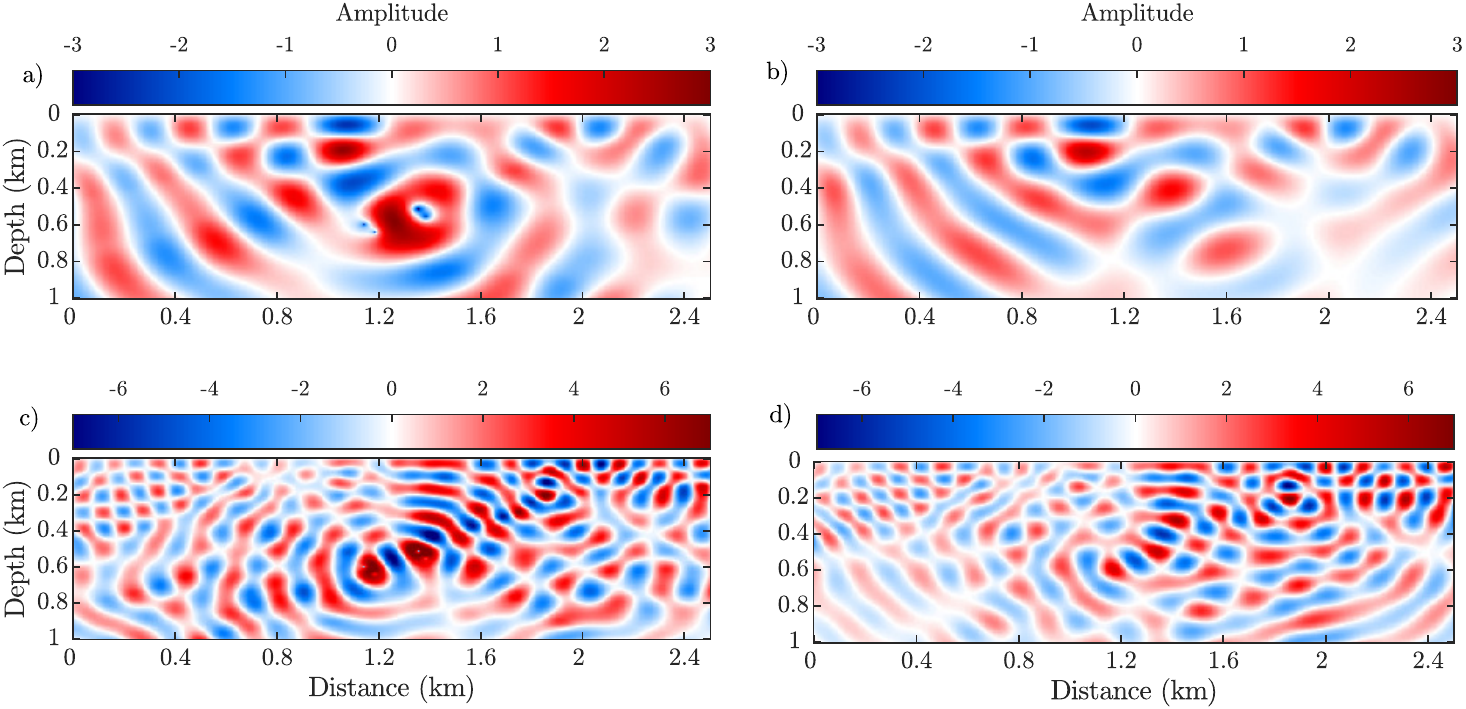}
\caption{Showing the ability of the data-assimilated wavefield to mimic the true wavefields when no information on the source is available. The Marmousi model with a source containing MSEs 1-4 is used. (a) True 7~Hz wavefield. (b) 7~Hz data-assimilated wavefield (Eq. \ref{close_U0}) with initial model 2. (c-d) Same as (a-b), but for 15~Hz.}
\label{fig_wavefield}
\end{figure}
\begin{figure}[htb]
\centering
\includegraphics[width=1\textwidth]{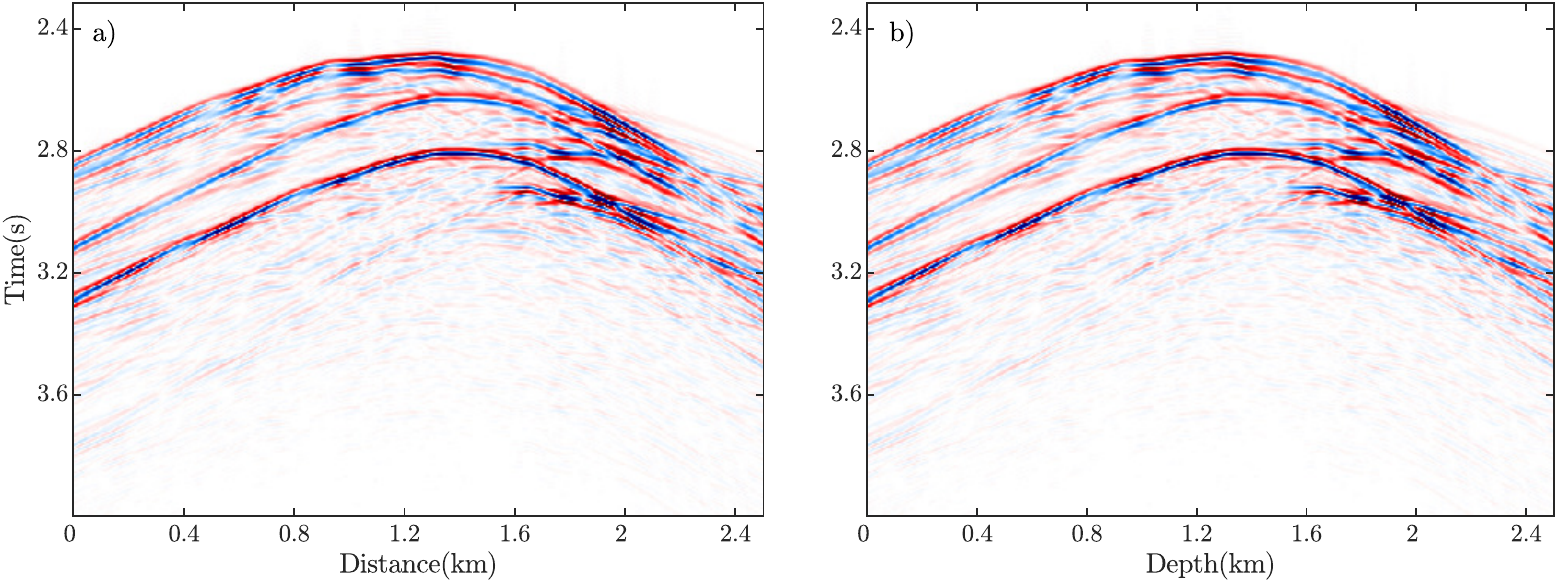}
\caption{Showing the match between the true and data assimilated seismograms when no information on the source is available. The Marmousi model with a source contains MSEs 1-4 is used. (a) True seismograms. (b) Seismograms created using data-assimilated wavefields (Eq. \ref{close_U0}) with initial model 2.}
\label{fig_shot_gather}
\end{figure}
%
\begin{figure}[htb]
\centering
\includegraphics[width=1\textwidth]{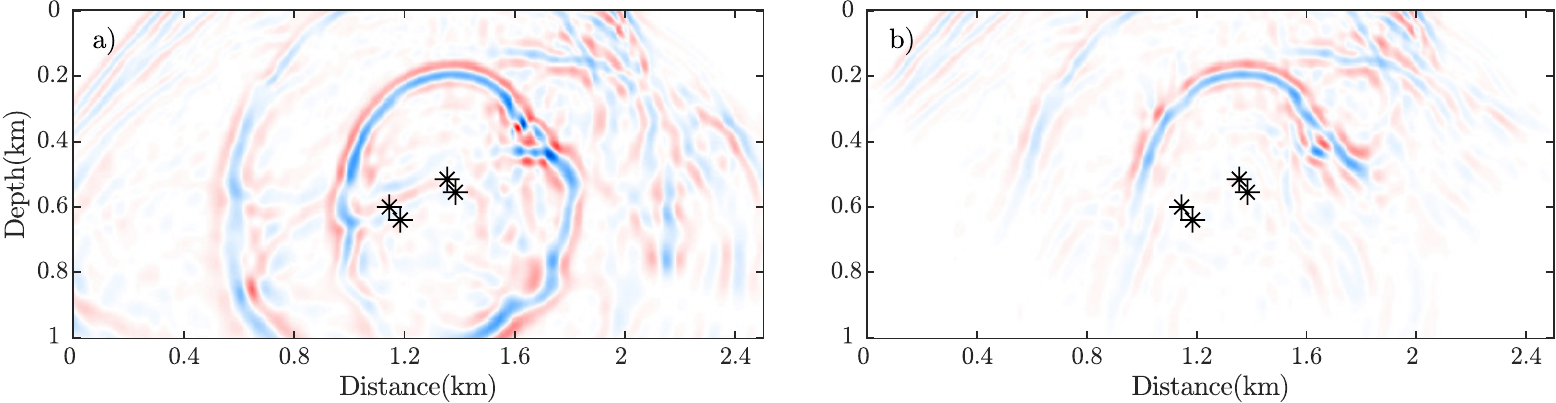}
\caption{Showing the ability of data-assimilated wavefield to capture many features of the true wavefield above the MSEs when no information on the source is available. The Marmousi model with a source contains MSEs 1-4 is used. A time slices at 2.7 s for (a) the true wavefield, (b) the data-assimilated wavefield with initial model 2.}
\label{fig_time_slice}
\end{figure}
We assess the method on the selected target of the synthetic Marmousi II model \citep{Martin_2006_M2E} with size $2.5~km \times 1~km$ when the grid spacing is $5~m$ (Fig. \ref{fig0}a). We use two different initial models for the tests: an accurate initial velocity model obtained by slightly smoothing the true model, referred to as initial model 1 (Fig. \ref{fig0}b), and a highly-smoothed, albeit kinematically accurate, initial velocity model referred to as initial model 2 (Fig. \ref{fig0}c). A direct comparison between the true model and initial models 1 and 2 at $X=1250~m$ is shown in Fig. \ref{fig0}d.
The different MSEs, which are used in this section, are shown in Fig. \ref{fig0}a using black asterisks and numbers. The MSEs of the right cluster are located at (520,1360)~m  [MSE 1] and (560,1380)~m [MSE 2], and those of the left cluster are located at (605,1150)~m [MSE 3], and (645,1175)~m [MSE 4]. The source signatures for MSEs 1 to 4 are Ricker wavelets with central frequencies [25, 31, 23, 29] Hz and central times [2.4, 2.56, 2.25, 2.2] s, respectively.\\   

\begin{figure}[htb]
\centering
\includegraphics[width=1\textwidth]{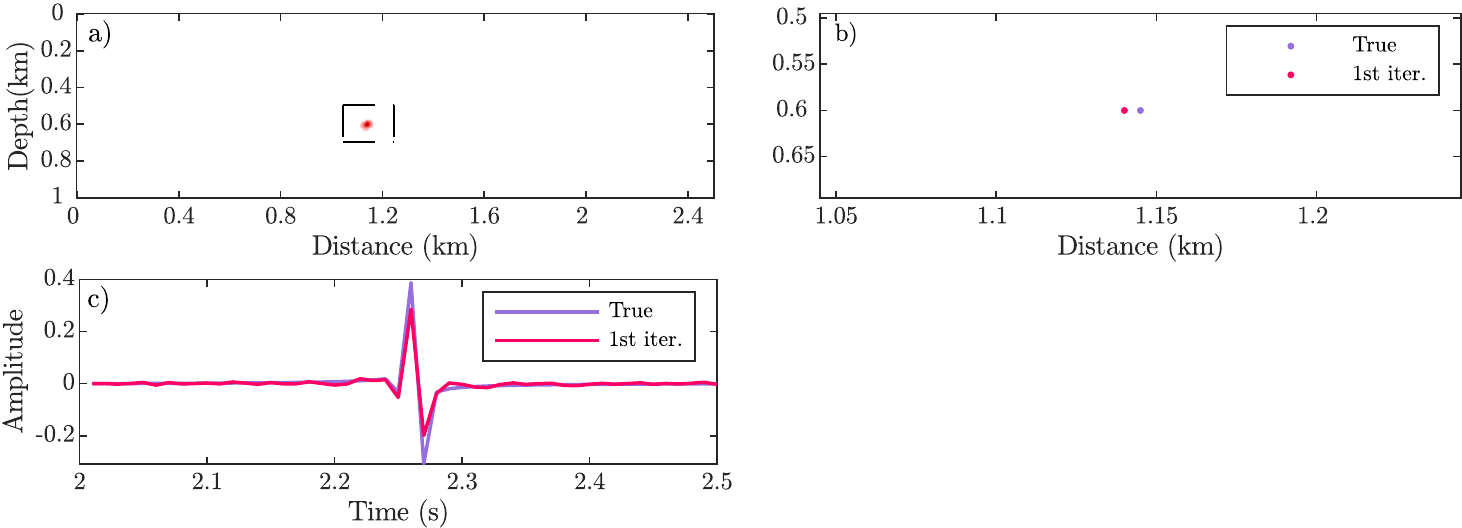}
\caption{Marmousi test with initial model 1 when the source just contains MSE 1. (a) The predicted source at the first outer iteration ($\bar{\bold{b}}^1$). (b) The zoomed result of peak finder applied on $\bar{\bold{b}}^1$. The margins of the zoomed part are shown using black dashed lines in (a). Also, the location of the true source is indicated by a purple point. (c) The true source signature (purple) and the estimated one (pink).}
\label{fig1}
\end{figure}

For all of the numerical tests, we compute wavefields with a nine-point stencil finite-difference method implemented with anti-lumped mass, where the stencil coefficients are optimized for each frequency \citep{Chen_2013_OFD}.  We use absorbing boundary conditions along the bottom, right and left sides of the model and a free-surface boundary condition at the surface.\\
We first illustrate how the initial data-assimilated wavefield is reconstructed (when no information on the source is available, eq. \ref{close_U0}). The source contains MSEs 1 to 4 and we use model 2 as the starting velocity model (Fig. \ref{fig0}c).  
Fig. \ref{fig_wavefield} shows the true frequency-domain wavefields and the initial data-assimilated wavefields for the 7~Hz and 15~Hz frequencies which are reconstructed by back-propagating the data, eq. \ref{close_U0}.
The ability of the data-assimilated wavefield to mimic the true wavefield can be assessed by comparing the left and right columns of this figure. This emphasizes how suitable the (transmission) microseismic configuration is to perform accurate wavefield reconstruction between the MSEs and the receivers using data assimilation in the initial model. The time-domain seismograms computed from the MSEs in the true model are shown in Fig. \ref{fig_shot_gather}a, while those computed in the initial model with data assimilation, when no information about the source is available, are shown in Fig. \ref{fig_shot_gather}b. These seismograms match almost perfectly because the data assimilated wavefield fits the data at the expense of the accuracy with which they satisfy the wave equation through the feedback term to the data (Eq. \ref{close_U0}). Also, a comparison between the time slices at 2.7 s of the true (Fig. \ref{fig_time_slice}a) and data-assimilated wavefields (Fig. \ref{fig_time_slice}b) shows that the data-assimilated wavefield captures many features of the true wavefield above the MSEs which again emphasizes the ability of data-assimilation to reconstruct an accurate wavefield. \\
We continue by assessing the performance of Algorithm \ref{Alg2cont0} in estimating the MSE location and signature estimation. We start by a simple setup and complicate it step by step. 
We perform the first test with a source that contains only MSE 1 using model 1 as the starting model (Fig. \ref{fig0}b). The source signature of this MSE is shown in Fig. \ref{fig1}c in purple. 
Note that all the wavelets in this section are created with the inverse Fourier transform of a few discrete frequencies, those that are used for inversion in IR-WRI. \\
We invert frequency components between 5~Hz and 45~Hz with a 2-Hz frequency interval simultaneously. We perform the inner loop for source location estimation (lines 5-8 of Algorithm \ref{Alg2cont0}) with 10 inner iterations ($n_l=9$).  The estimated $\bar{\bold{b}}^1$ at the final iteration of the inner loop is shown in Fig. \ref{fig1}a. 
%
%
The event location found by the peak finder algorithm is shown in Fig. \ref{fig1}b in pink while the true one is shown in purple. The results are zoomed in from the full model; the limits of the zoom are shown with the black dashed lines in Fig. \ref{fig1}a. We observe that the estimated MSE location is close to the true one because of the accuracy of the initial velocity model. Finally, the estimated source signature is shown in Fig. \ref{fig1}c (pink curve), which matches the true one well. \\

%
%
%
\begin{figure}[htb]
\centering
\includegraphics[width=1\textwidth]{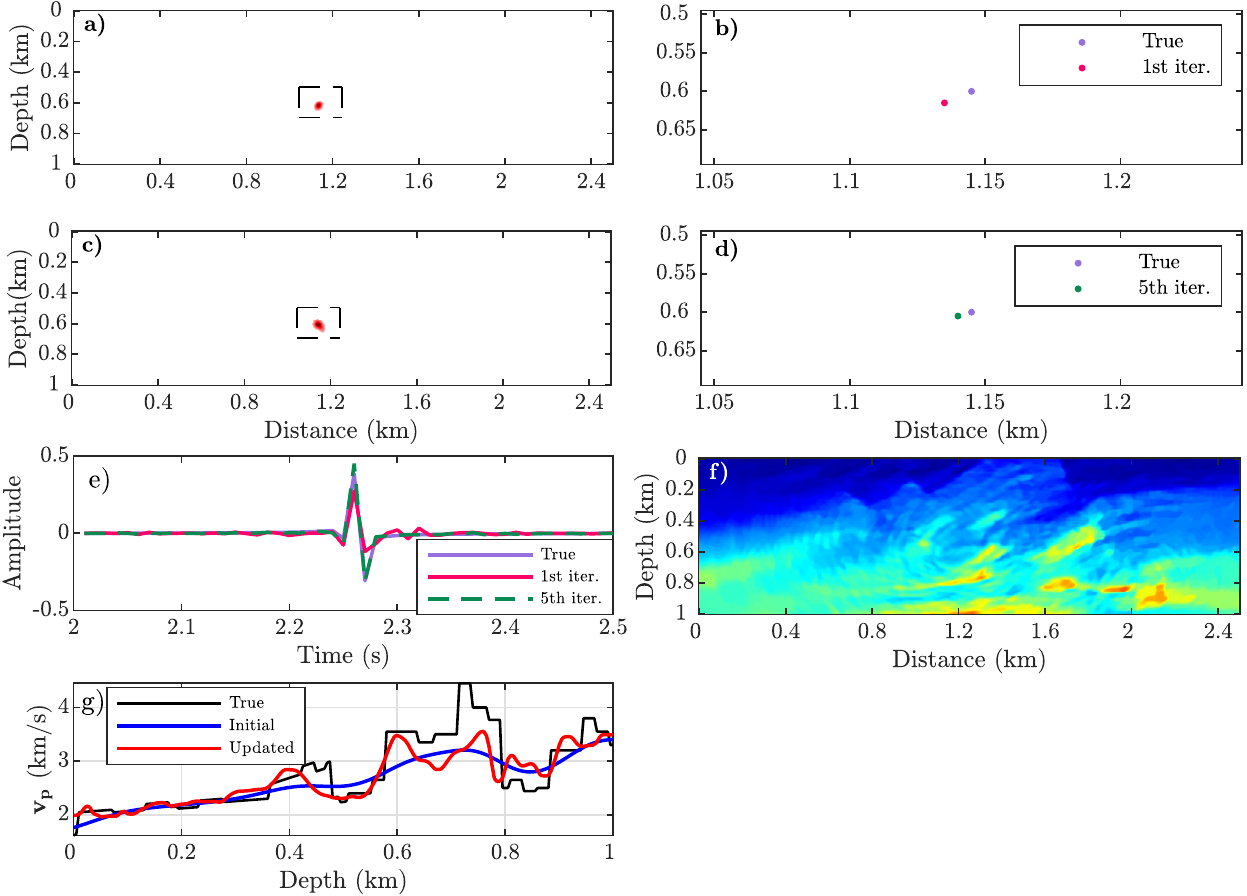}
\caption{Marmousi test with initial model 2 when the source just contains MSE 1. (a) Predicted source at the final inner iteration of outer iteration 1 ($\bar{\bold{b}}^1$). (b) The zoomed result of peak finder applied on $\bar{\bold{b}}^1$. The margins of the zoomed part are shown using black dashed lines in (a). Also, the location of the true source is indicated by a purple point. (c-d) Same as (a-b) but for $\bar{\bold{b}}^5$.  (e) The true source signature (purple) and the estimated at outer iteration 1 (pink) and iteration 5 (green). (f) Updated velocity model after five iterations.}
\label{fig2}
\end{figure}

We make this test more complicated by starting with initial velocity model 2. The estimated predicted source at the final inner iteration of outer iteration 1 ($\bar{\bold{b}}^1$) is shown in Fig. \ref{fig2}a and the result of the application of the peak finder algorithm is shown in Fig. \ref{fig2}b in pink. We see that the error in the estimated location is more than the previous case (Fig. \ref{fig1}b) where the initial model was accurate. Also, the estimated source signature is shown in pink in Fig. \ref{fig2}e. We can clearly see that neither the estimated MSE location nor the estimated signature are accurate. To improve the quality of these results, we should update the velocity model. The estimated predicted source after five outer iterations (at the final iteration of inner loop), $\bar{\bold{b}}^5$, and the result of peak finder algorithm are shown in Figs. \ref{fig2}c-\ref{fig2}d, respectively. Finally, the estimated source signature at outer iteration 5 is shown in green in Fig. \ref{fig2}e. It can be seen that the quality of the estimated location and signature improved significantly. The updated velocity model at outer iteration five is shown in Fig. \ref{fig2}f. Although the velocity update is not useful for geological interpretation due to the limited illumination provided by the single event, it improves the structures between the event and the receivers and, as a result, leads to an improved MSE location and signature estimation. \\

%
%
\begin{figure}[htb]
\centering
\includegraphics[width=1\textwidth]{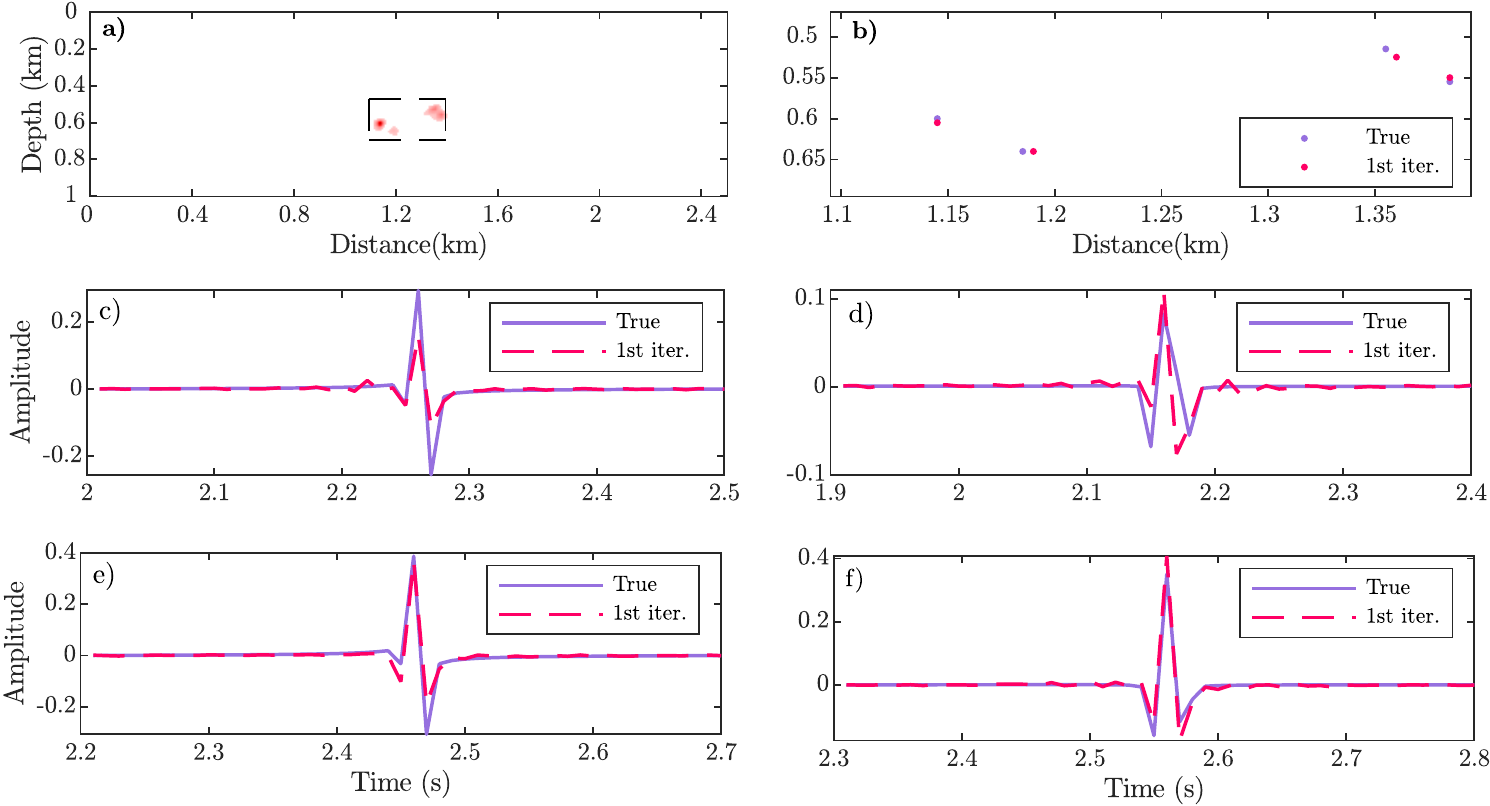}
\caption{The Marmousi test with initial model 1 when the source contains MSEs 1-4. (a-b) $\bar{\bold{b}}^1$ and its selected peaks using the peak finder. (c-f) True source signatures (purple) and the estimated ones (pink) for MSEs 1 to 4, respectively.}
\label{fig3}
\end{figure}

\begin{figure}[htb]
\centering
\includegraphics[width=1\textwidth]{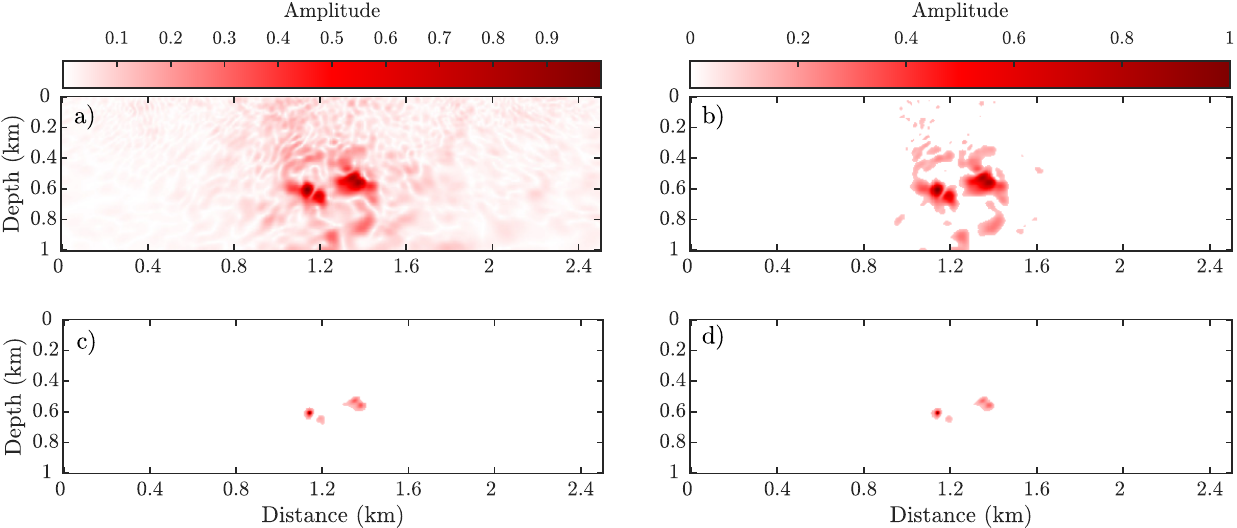}
\caption{Showing the impact of sparsifying regularization on the estimated predicted sources for the Marmousi test with initial model 1 when the source contains MSEs 1-4. (a-b) The predicted source at $k=0~,l=0$ (a) before, (b) after Berhu regularization. (c-d) Same as (a-b) but for $k=0~,l=9$. $k$ and $l$ refer to the outer and inner iteration of algorithm~\ref{Alg2cont0}, respectively.}
\label{fig3_1}
\end{figure}

To make the test more representative of the microseismic scenario, we repeat the Marmousi test when the source gathers MSEs 1-4 (Fig. \ref{fig0}a). We start with the initial model 1. The predicted source after 1 outer iteration ($\bar{\bold{b}}^1$) and the result of peak finder are shown in Figs. \ref{fig3}a-\ref{fig3}b, respectively. We see that, because of the good initial model, there are just 4 peaks of energy in Fig. \ref{fig3}. Also, the estimated source signatures for MSEs 1 to 4 are shown in Figs. \ref{fig3}c-\ref{fig3}f, respectively, in pink. We have accurately recovered both the locations and signatures for all MSEs. \\
%
%
%
%
%
Even with a good initial model, the estimated predicted source has a lot of unwanted energy, and sparsifying regularization is necessary to find the locations of the MSEs from the predicted source. The predicted source for this test at inner iteration 1 ($k=0,~l=0$) before and after applying Berhu regularization are depicted in Figs. \ref{fig3_1}a-\ref{fig3_1}b, respectively. We can see the improvement, although it is not enough to have an accurate picking for the location of the MSEs. The same results after 10 inner iterations ($k=0,~l=9$) are shown in Figs. \ref{fig3_1}c-\ref{fig3_1}d. First, we see that the regularization has little effect on the predicted source in Fig. \ref{fig3_1}d  (compare Fig. \ref{fig3_1}c) because it is already sparse. Second, we see the improvement in Fig. \ref{fig3_1}d compared to Fig. \ref{fig3_1}b, where the regularization has significantly improved the MSE location selection.\\

%
%
%
%
%
%
\begin{figure}[htb]
\centering
\includegraphics[width=1\textwidth]{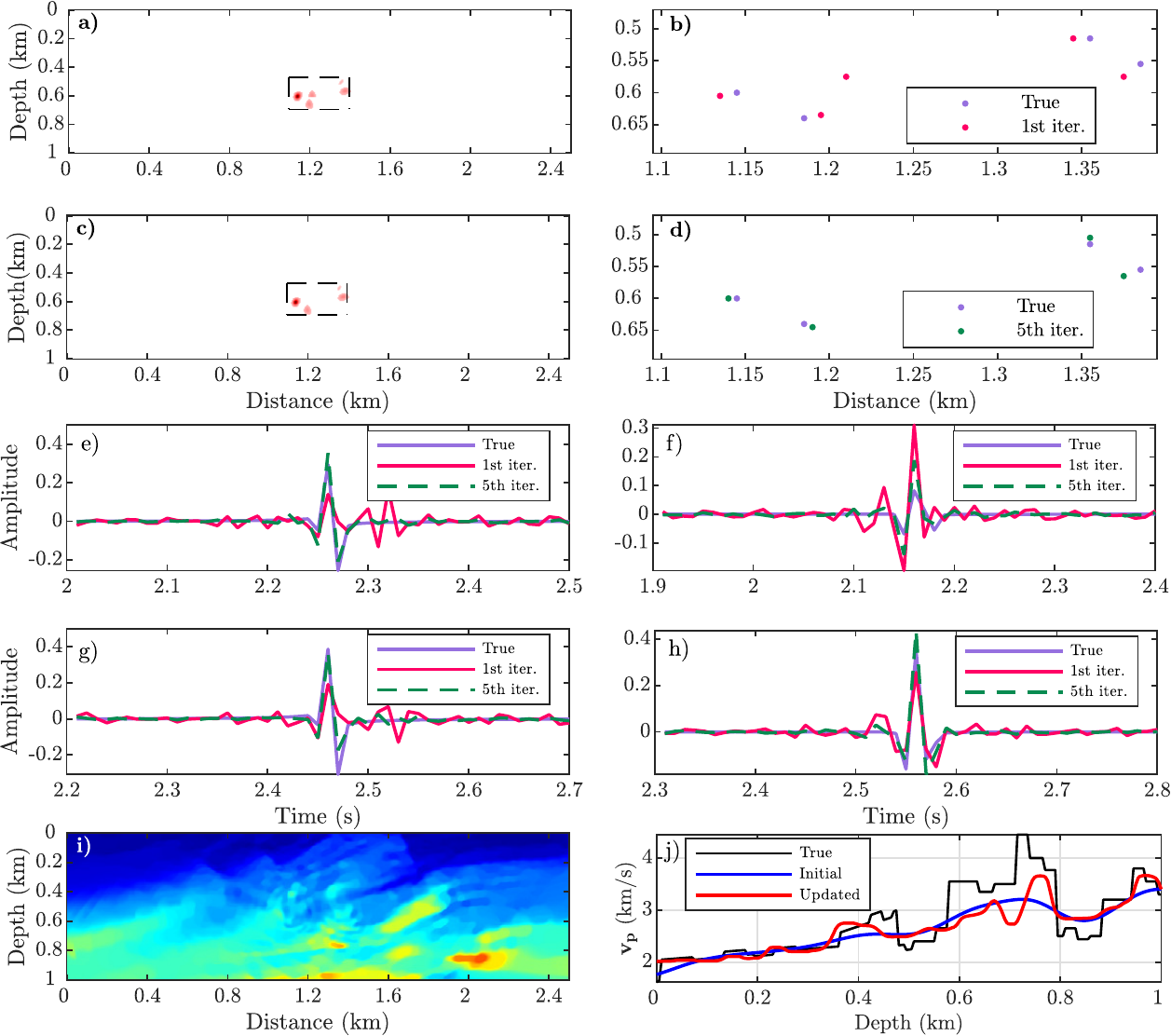}
\caption{The Marmousi test with initial model 2 when the source contains MSEs 1-4. (a-b) $\bar{\bold{b}}^1$ and its selected peaks using peak finder. (c-d) $\bar{\bold{b}}^5$ and its selected peaks. (e-h) True MSE signatures (purple) and the estimated ones for MSEs 1 to 4, respectively, at outer iteration one (pink) and iteration five (green). (i) The updated velocity model after five iterations.}
\label{fig4}
\end{figure}
\begin{figure}[htb]
\centering
\includegraphics[width=1\textwidth]{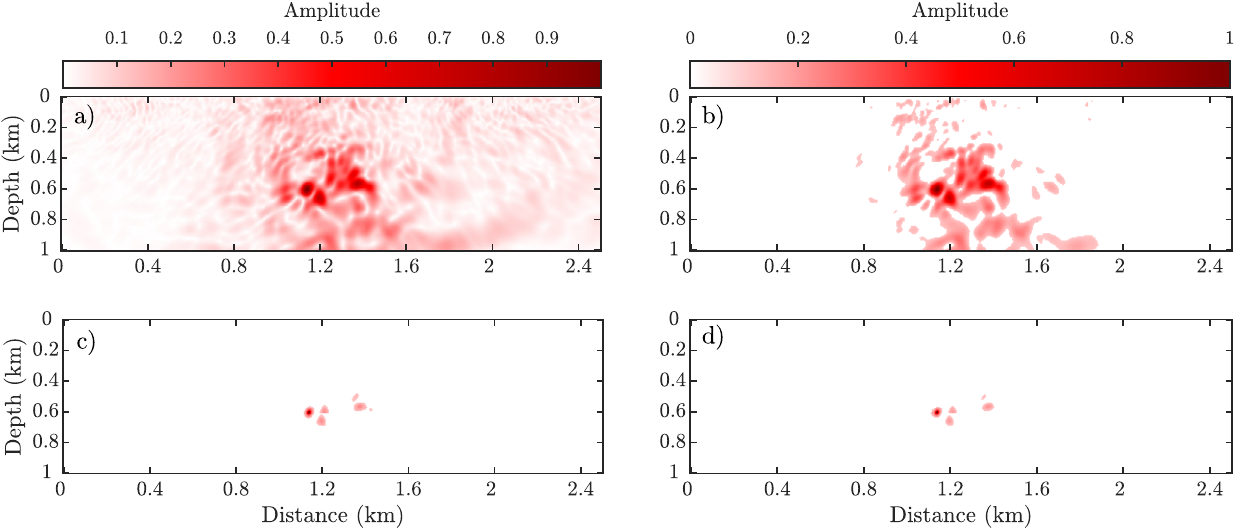}
\caption{Same as Fig. \ref{fig3_1}, but for Marmousi test with initial model 2 when the source contains MSEs 1-4.}
\label{fig4_1}
\end{figure}

We repeat the Marmousi test with initial velocity model 2 when the source contains MSEs 1-4, and we use five outer iterations to update the velocity model to improve the quality of the estimated MSEs. The estimated $\bar{\bold{b}}^1$ (at the final inner iteration) and its located peaks are shown in Figs. \ref{fig4}a-\ref{fig4}b and the same results for outer iteration five are shown in Figs. \ref{fig4}c-\ref{fig4}d. At outer iteration 1, the peak finder algorithm finds more peaks (five) than the true number of MSEs because of the inaccurate velocity model (Fig. \ref{fig4}b). Then, the data-assimilated wavefield and source signatures for all these selected MSEs are updated jointly using Eq. \ref{closed_u_2}.  Because the selected MSE at the middle of Fig. \ref{fig4}b does not contribute to the recorded data, the algorithm finds a small source signature for this fake MSE. As soon as the velocity model improves, we have a better predicted source (Fig. \ref{fig4}c), and the peak finder algorithm selects four points close to the true MSEs (Fig. \ref{fig4}d).  Also, the estimated MSE signatures at outer iterations 1 and 5 are shown in Figs. \ref{fig4}e-\ref{fig4}h, respectively, for MSEs 1 to 4 in pink for outer iteration one and in green for outer iteration five. Finally, the updated velocity model after five outer iterations is shown in Fig. \ref{fig4}i.  First, the quality of the estimated signatures of the MSEs are improved by updating the velocity model. Second, the updated velocity model captures the trends of dominant structures, and is improved compared to the initial model 2.  \\
Like the previous test, the estimated source locations at inner iteration 1 ($k=0,~l=0$) before and after applying Berhu regularization are depicted in Figs. \ref{fig4_1}a-\ref{fig4_1}b, respectively, and for inner iteration 10 ($k=0,~l=9$) shown in Figs. \ref{fig4_1}c-\ref{fig4_1}d. We see that the energy is less focused compared to the case with an accurate initial model (Figs. \ref{fig3_1}a-\ref{fig3_1}b). We also see that the regularization significantly improves the predicted source map, making it ready for MSEs location picking.\\
%
%
%
%
%
%
%
%
%
%
%
%
\begin{figure}[htb]
\centering
\includegraphics[width=1\textwidth]{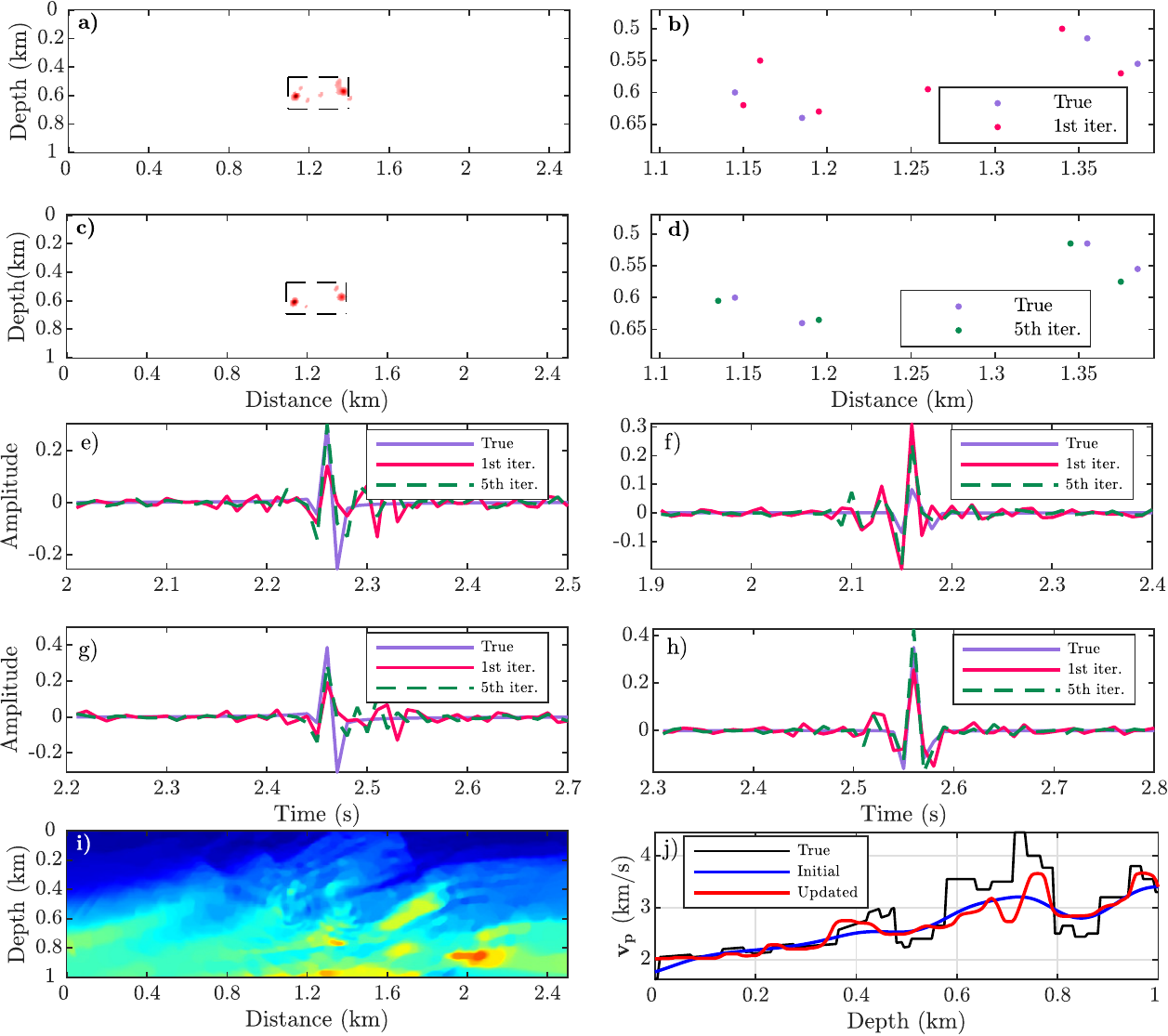}
\caption{Same as Fig. \ref{fig4} but with noisy data with SNR=$5~db$.}
\label{fig4_noisy}
\end{figure}

We assess now the robustness of the method against noise. We repeat the previous test (with initial model 2 and the source that contains MSEs 1-4) when Gaussian distributed random noise with SNR=$5~db$ is added to the data. We use the same configuration as for the noiseless test (Fig. \ref{fig4}). The results are shown in Fig. \ref{fig4_noisy}. In comparison with the noiseless case, we see that the predicted source has more peaks than the true number of MSEs at the first outer iteration (Fig. \ref{fig4_noisy}a), but it is improved at the fifth outer iteration (Fig. \ref{fig4_noisy}c). Also, we see that the estimated MSEs are not significantly changed in Fig. \ref{fig4_noisy}c compared to the noiseless case (Fig. \ref{fig4}c). For the estimated signatures, Figs. \ref{fig4_noisy}e-\ref{fig4_noisy}h, we see the quality of the estimated signatures are degraded compare to the noiseless case (Figs. \ref{fig4}e-\ref{fig4}h). Finally, the updated velocity model for this noisy test, Fig. \ref{fig4_noisy}i, is close to the noiseless one in Fig. \ref{fig4}i.\\

\section{Discussion}
%
%
%
%
%
%
%
%
%
%
%

%

In this paper, we have proposed an algorithm based on ADMM-based FWI for characterizing weak seismic events that is valid at any scale. We have focused on microseismic imaging to find the location and signature of MSEs, but the method is general and could be applied to tectonic events or any other situation in which the source location is unknown. The proposed algorithm does not require any information or assumptions about the sources, and it is able to find the location and time signature of seismic events and update the background velocity model if required. The proposed algorithm consists of three steps: 1- finding the number and location of seismic events, 2- jointly updating the data-assimilated wavefield and the signature of seismic events, 3- updating the background velocity model provided that the recorded data provides sufficient illumination of the model.\\

If the location of the seismic events are known, \citet{Fang_2018_SEF,Aghamiry_2021_EES} have shown that it is possible to jointly estimate the data-assimilated wavefield and the event signatures with high accuracy by solving a linear least-squares problem. But, when the location and the number of seismic events are unknown, the first step of the proposed algorithm tries to find this information by applying a peak finder algorithm on the sparsified predicted source map. 
By using Eq. \eqref{close_U0} and considering the estimated source as $\A(\m^0,\w)\u(\w)^{0}$, the predicted source locations are generated by first propagating the data backward in time from the receiver positions and then the blurring effects induced by the limited bandwidth of the data and the limited spread of the receivers are corrected with the sparsity-promoting regularization. The average of these estimated sources helps to find the number and the location of seismic events. However, the success of this step strongly depends on the sparsifying regularization that is applied on the estimated sources.\\ 
At the second step of the algorithm, the seismic event signatures and data-assimilated wavefield are jointly updated from the estimated locations. When the algorithm finds incorrect seismic events in the first step, one contribution of this second step is to mitigate their footprint by assigning low-amplitude signatures to them.
This wavefield and signature refinement, as well as the updating of the background model, help the algorithm to better-estimate the source and seismic event locations in subsequent outer iterations.\\

The numerical results show that we should have at least a kinematically accurate initial model to have a good estimate of the seismic events. But, if the model is kinematically incorrect, we should have sufficient illumination of the model in the recorded data to update it during the outer iteration.  We note that we cannot update the velocity model everywhere, but the region in which we can update the model is the region through which the recorded waves travel.  Thus we can update the part of the model that is most important for the estimation of the source. \\

The proposed algorithm is in the frequency domain when the inversion is limited to a few frequencies. The computational burden of the proposed algorithm is primarily in solving the data-assimilated subproblem with a direct solver, Eqs. \ref{close_UK}, \ref{closed_u_2}. For problems with large computational domains, these subproblems can be solved efficiently with preconditioned iterative solvers. \citet{Rezaei_2021_ALB} show that conjugate-gradient with additive Schwarz domain-decomposition preconditioner \citep{Dryja_1987_AAV} has the best performance for solving data-assimilated subproblem in comparing it to other solvers. The computational burden of all the other steps of the algorithm are negligible. 
%
%
\section{CONCLUSIONS} 
We proposed a method based on ADMM-based FWI for characterizing weak seismic events at different scales. 
When the source is added to the unknowns of FWI, in addition to the wavefields and model parameters, it becomes a severely underdetermined problem, and it is challenging to uniquely determine the source without prior information. 
We present a method to solve this problem that does not require prior information about the sources, although it does require a kinematically correct velocity model.  The proposed method consists of three steps:
With a finite band of frequencies and appropriate sparsifying regularizations on the source, seismic events are selected with a peak finder algorithm applied on the estimated average (over frequency) source.
The time signatures and wavefields are jointly updated for the selected seismic events by solving a linear least-squares problem.
The velocity model is updated by applying appropriate regularization on the model.  
We validate the proposed method with synthetic tests for microseismic event characterizations as a proof of concept.  The method works well, even for closely spaced events and has acceptable performance with extremely low signal-to-noise ratio data. The method can be tailored to earthquake relocation or microseismic monitoring. Ongoing work involves extending the method to 3D elastic physics and moment tensor estimation in the prospect of real data applications.
\section*{ACKNOWLEDGMENTS}  
This study was partially funded by the WIND consortium (\textit{https://www.geoazur.fr/WIND}), sponsored by Chevron, Shell and Total. The authors are grateful to the OPAL infrastructure from 
Observatoire de la Côte d'Azur (CRIMSON) for providing resources and support. This work was granted access to the HPC resources of IDRIS under the allocation A0050410596 made by GENCI.
A. Malcolm acknowledges the NSERC Discovery Grant Program as well as Chevron, InnvoateNL and the NSERC Industrial Research Chair program.

\bibliographystyle{gji}
\newcommand{\SortNoop}[1]{}

\end{document}